\documentclass[12pt]{article}
\usepackage{amsmath}
\usepackage{latexsym}
\usepackage{amssymb}
\usepackage{a4wide}
\newtheorem{thm}{Theorem}[section]
\newtheorem{la}[thm]{Lemma}
\newtheorem{Defn}[thm]{Definition}
\newtheorem{Remark}[thm]{Remark}
\newtheorem{Note}[thm]{Note}
\newtheorem{prop}[thm]{Proposition}

\newtheorem{Example}[thm]{Example}
\newtheorem{Examples}[thm]{Examples}
\newtheorem{Problems}[thm]{Problems}

\newtheorem{Problem}[thm]{Problem}
\newtheorem{Number}[thm]{\!\!}
\newenvironment{defn}{\begin{Defn}\rm}{\end{Defn}}

\newenvironment{example}{\begin{Example}\rm}{\end{Example}}

\newenvironment{rem}{\begin{Remark}\rm}{\end{Remark}}

\newenvironment{proof}{{\noindent\bf Proof.}}%
                  {\nopagebreak\hspace*{\fill}$\Box$\medskip\medskip\par}   
\newcommand{\Punkt}{\nopagebreak\hspace*{\fill}$\Box$}

\newcommand{\ve}{\varepsilon}
\newcommand{\at}{\symbol{'100}}

\newcommand{\wt}{\widetilde}

\newcommand{\n}{\rm}
\newcommand{\impl}{\Rightarrow}

\newcommand{\mto}{\mapsto}

\newcommand{\N}{{\mathbb N}}
\newcommand{\R}{{\mathbb R}}

\newcommand{\Z}{{\mathbb Z}}
\newcommand{\C}{{\mathbb C}}
\newcommand{\K}{{\mathbb K}}

\newcommand{\cV}{{\cal V}}

\newcommand{\Diff}{{\mbox{{\rm Diff}}}}
\newcommand{\dl}{{\displaystyle \lim_{\longrightarrow}}}
\newcommand{\pl}{{\displaystyle \lim_{\longleftarrow}}}

\newcommand{\take}{\backslash}
\newcommand{\sub}{\subseteq}

\newcommand{\im}{\mbox{\n im}}

\newcommand{\pr}{\mbox{\rm pr}}

\newcommand{\id}{\mbox{\n id}}

\newcommand{\Supp}{\mbox{\n supp}}
\newcommand{\ev}{\mbox{\n ev}}
\newcommand{\comp}{{\mbox{\n \footnotesize comp}}}
\newcommand{\loc}{{\mbox{\n \footnotesize loc}}}

\begin{document}
\begin{center}
{{\Large\bf Discontinuous Non-Linear Mappings on\vspace{1mm}
Locally Convex Direct Limits}}\\[2.5mm]
Helge Gl\"{o}ckner
\end{center}
\begin{abstract}
\hspace*{-7.2 mm}
We show that the self-map $f\!: C^\infty_c(\R)\to C^\infty_c(\R)$,
$f(\gamma):=\gamma\circ \gamma\, -\, \gamma(0)$
of the space of real-valued test functions on the line
is discontinuous, although its restriction
to the space $C^\infty_K(\R)$ of functions
supported in~$K$ is smooth (and hence continuous), for each
compact subset~$K\sub \R$.
More generally, we construct mappings with analogous
pathological properties on spaces of compactly supported
smooth sections in vector bundles over non-compact bases.
The results are useful in infinite-dimensional Lie theory,
where they can be used to analyze the precise direct limit
properties of test function groups and groups of
compactly supported diffeomorphisms.
%
%
%
%
%
\end{abstract}
\vspace{2mm}
\begin{center}
{\bf\Large Introduction}
\end{center}
Let $E_1\sub E_2\sub \cdots$ be an ascending sequence of
locally convex spaces which does not become
stationary,
and such that
$E_{n+1}$ induces the given topology
on~$E_n$, for each~$n$.
It is a well-known phenomenon
that the topology on $E:=\bigcup_{n\in\N}E_n$
making~$E$ the direct limit of the spaces~$E_n$
in the category of locally convex spaces (and continuous
linear maps) can be properly coarser
than the topology making~$E$
the direct limit of its subspaces~$E_n$
in the category of topological spaces.
For example, this phenomenon
occurs whenever
each~$E_n$ is an infinite-dimensional
Fr\'{e}chet space (cf.\ \cite[Prop.\,4.26\,(ii)]{KaM}).
In particular, the locally convex direct
limit topology on the space $C^\infty_c(\R)=
\dl\, C^\infty_{[-n,n]}(\R)$\vspace{-1.3 mm}
of test functions
is properly coarser than the topology of direct
limit topological space
(cf.\ also \cite[p.\,506]{Dud}).\\[3.5mm]
So, for abstract reasons,
discontinuous mappings
on the space of test functions $C^\infty_c(\R)$
are known to exist
whose restriction to $C^\infty_{[-n,n]}(\R)$
is continuous for each $n\in \N$.
In this article, we describe
such a mapping explicitly,
whose restriction to $C^\infty_{[-n,n]}(\R)$
is not only continuous but actually smooth
(Proposition~\ref{prototype}).
More generally, for
every $\sigma$-compact,
non-compact, finite-dimensional smooth manifold~$M$
of positive dimension
and locally convex space $E\not=\{0\}$,
we construct
a discontinuous map
$f\!:C^\infty_c(M,E)\to C^\infty_c(M,\R)$
whose
restriction
to $C^\infty_K(M,E)$ is smooth,
for each compact subset~$K$ of~$M$.
An analogous result is obtained
for the space~$C^\infty_c(M,E)$ of compactly supported smooth sections
in a bundle of locally convex spaces
$E\to M$ over~$M$, with non-trivial fibre
(Theorem~\ref{fingeruebung}).\\[3.5mm]
{\bf Further developments.}
The preceding result is useful for
the investigation
of direct limit properties
of infinite-dimensional
Lie groups.
As shown in~\cite{DL3},
it entails that
there are discontinuous
(and hence non-smooth)
mappings on the Lie group $\Diff_c(M)
=\bigcup_K \text{Diff}_K(M)$
of compactly supported smooth diffeomorphisms
of~$M$ (as in \cite{Mic} or \cite{DIF}),
whose restriction to $\Diff_K(M):=\{\phi\in \Diff(M)\!:
\phi|_{M\setminus K}=\id_{M\setminus K}\}$
is smooth, for each compact subset $K\sub M$.
A similar pathology occurs for the
Lie group $C^\infty_c(M,G)=\bigcup_K C^\infty_K(M,G)$
of compactly supported smooth maps
with values in a non-discrete
finite-dimensional Lie group (as in~\cite{GCX}).
In this way, we obtain one half of the following
table, which describes whether
$\text{Diff}_c(M)=\dl\, \text{Diff}_K(M)$\vspace{-.5mm} and
$C^\infty_c(M,G)=\dl\, C^\infty_K(M,G)$\vspace{-.5mm}
holds in the categories shown:\vspace{2mm}
\begin{center}
\begin{tabular}{||c||c|c||}\hline\hline
category $\backslash$ group &
$C^\infty_c(M,G)$ & $\text{Diff}_c(M)$ \\ \hline\hline
Lie groups & yes & yes \\ \hline
topological groups & yes & yes \\ \hline
smooth manifolds & no & no \\ \hline
topological spaces & no & no \\ \hline\hline
\end{tabular}\vspace{5mm}
\end{center}
For the proof, see~\cite{DL3}
(cf.\ also \cite{TSH} for related results).\\[3.5mm]
The present constructions of pathological
mappings are complemented
by investigations in
\cite{ZOO}--\cite{DIF}
(cf.\ also \cite{DRn}).
In these articles,
a mild additional
property is introduced
which ensures
that a map
$f\!: C^\infty_c(M,E)\to C^\infty_c(N,F)$
between spaces of test functions
(or compactly supported sections)
satisfying this property
(an ``almost local'' map)
is indeed smooth if and only
if it is smooth on $C^\infty_K(M,E)$
for each~$K$.
In contrast to these mappings,
the pathological examples presented here are
extremely non-local.\\[3mm]
In the final section, we describe
examples of discontinuous bilinear
mappings which are continuous (and hence analytic)
on each step of a directed sequence of subspaces.
\section{Preliminaries}\label{secprelim}
In this article,
we are working in the setting
of infinite-dimensional differential
calculus known as Keller's $C^\infty_c$-theory,
based on smooth maps in the sense
of Michal-Bastiani
(see \cite{BED},
\cite{Ham}, \cite{Mic}, \cite{Mil}
for further information).
\begin{defn}
Let $E$, $F$ be locally convex spaces and $f\!:U\to F$
be a mapping, defined on an open subset~$U$ of~$E$.
We say that $f$ is {\em of class~$C^0$\/}
if~$f$ is continuous.
If~$f$ is a continuous map such that
the two-sided directional derivatives
\[
df(x,v)=\lim_{t\to 0} {\textstyle
\frac{1}{t}\left( f(x+tv)-f(x)\right)}
\]
exist for all $(x,v)\in U\times E$,
and the map $df\!: U\times E\to F$
so defined is continuous,
then $f$ is said to be
{\em of class~$C^1$\/}. Recursively, given $k\in \N$
we call~$f$ a mapping of class~$C^{k+1}$
if it is of class~$C^1$ and $df$ is of class~$C^k$
on the open subset $U\times E$
of $E\times E$.
We set $d^{k+1}f:=d(d^kf)=d^k(df)\!:
U\times E^{2^{k+1}-1}\to F$
in this case. The function~$f$ is called
{\em smooth\/} (or of class $C^\infty$)
if it is of class~$C^k$ for each $k\in \N_0$.
\end{defn}
\begin{defn}
Let $M$ be a finite-dimensional, $\sigma$-compact
smooth manifold and~$E$ be a locally convex topological vector space.
We equip the vector space $C^\infty(M,E)$ of
$E$-valued smooth mappings~$\gamma$ on~$M$
with the topology of uniform convergence of
$\partial^\alpha(\gamma\circ \kappa^{-1})$
on compact subsets of~$V$,
for each chart $\kappa\!:M\supseteq U\to V\sub \R^d$
of~$M$ and multi-index $\alpha\in \N_0^d$
(where $d:=\dim(M)$).
Given a compact subset $K\sub M$, we equip the
vector subspace
$C^{\,\infty}_K(M,E):=\{\gamma\in C^\infty(M,E)\!:
\gamma|_{M\take K}=0\}$ of $C^\infty(M,E)$ with the induced
topology. We give
$C^\infty_c(M,E):=\bigcup_K C^{\, \infty}_K(M,E)=\dl\,
C^{\,\infty}_K(M,E)$\vspace{-.8 mm}
the locally convex direct limit topology.
We abbreviate
$C^\infty_c(M):=C^\infty_c(M,\R)$,
$C^\infty(M):=C^\infty(M,\R)$,
and $C^{\, \infty}_K(M):= C^{\,\infty}_K(M,\R)$.
Further details can be found,
e.g., in~\cite{GCX}.
\end{defn}
\section{Example of a discontinuous mapping on {\boldmath
$C^\infty_c(\R)$}}\label{secline}
We show that the map
$f\!: C^\infty_c(\R)\to C^\infty_c(\R)$,
$\gamma\mto \gamma\circ \gamma-\gamma(0)$
is discontinuous,
although its restriction to
$C^\infty_{[-n,n]}(\R)$ is smooth, for each $n\in \N$.\\[2mm]
The following
fact is essential
for our constructions.
It follows from \cite[Cor.\,3.13]{KaM}
and is also a special case
of~\cite[Prop.\,11.3]{ZOO}.
For the convenience of the reader,
we offer a direct, elementary proof
as an appendix.
\begin{la}\label{La1}
The composition map
\[
\Gamma\!: C^\infty(\R^n,\R^m)\times
C^\infty(M,\R^n)\to C^\infty(M,\R^m)\,,\qquad
\Gamma(\gamma,\eta)\, :=\, \gamma\circ \eta
\]
is smooth,
for each finite-dimensional, $\sigma$-compact
smooth manifold~$M$ and $m,n\in \N_0$.\Punkt
\end{la}
For the following proof,
recall that
the sets
\[
{\textstyle \cV(k,e):=\left\{
\gamma\in C^\infty_c(\R)\!:\;
(\forall n\in \Z)\,
(\forall j\in \{0,\ldots, k_n\})\,
(\forall x\in [n-\frac{1}{2},n+\frac{1}{2}])\;
|\gamma^{(j)}(x)|<\ve_n\right\}}
\]
form a basis of open
zero-neighbourhoods
for the topology on $C^\infty_c(\R)$,
where $k=(k_n)\in (\N_0)^\Z$ and $e=(\ve_n)\in (\R^+)^\Z$
(cf.\ \cite[\S\,II.1]{Sch};
see \cite[Prop.\,4.8]{GCX}).
\begin{prop}\label{prototype}
$f\!: C^\infty_c(\R)\to C^\infty_c(\R)$,
$\gamma\mto \gamma\circ \gamma-\gamma(0)$
has the following properties:
\begin{itemize}
\item[\n (a)]
The restriction of~$f$ to a map
$C^\infty_{[-n,n]}(\R)\to C^\infty_c(\R)$
is smooth $($and hence continuous$)$, for each $n\in \N$.
\item[\n (b)]
$f$ is discontinuous at $\gamma=0$.
\end{itemize}
\end{prop}
\begin{proof}
(a) Fix $n\in \N$;
we have to show
that $f|_{C^\infty_{[-n,n]}(\R)}
\!: C^\infty_{[-n,n]}(\R)\to C^\infty_c(\R)$
is smooth.
The image of this map being contained
in the closed vector subspace
$C^\infty_{[-n,n]}(\R)$
of $C^\infty_c(\R)$,
which also
is a closed vector subspace
of $C^\infty(\R)$ (with the same induced topology),
it suffices to show
that the map $C^\infty_{[-n,n]}(\R)\to
C^\infty(\R)$,
$\gamma\mto \gamma\circ \gamma-\gamma(0)$
is smooth
(see \cite[Prop.\,1.9]{SEC} or \cite[La.\,10.1]{BGN}).
Now $\gamma\mto \gamma(0)$
being a continuous linear (and thus smooth)
map, it suffices to show that
$C^\infty_{[-n,n]}(\R)\to C^\infty(\R)$,
$\gamma\mto \gamma\circ \gamma$ is smooth.
This readily follows from Lemma~\ref{La1}.

(b) Consider the zero-neighbourhood
$V:=\cV((|n|)_{n\in \Z}, (1)_{n\in \Z})$
in $C^\infty_c(\R)$.
Let $k=(k_n)\in (\N_0)^\Z$
and $e=(\ve_n)\in (\R^+)^\Z$ be arbitrary.
We show that $f(\cV(k,e))\not\sub V$.
Since $f(0)=0$,
this entails
that~$f$ is discontinuous
at $\gamma=0$.
It is easy
to construct a function
$h\in C^\infty_c(\R)$
such that $\Supp(h)\sub \;
]{- \frac{1}{2}},\frac{1}{2}[$
and $h(x)=x^{k_0+1}$
for all $x\in [-\frac{1}{4},\frac{1}{4}]$.
Then $rh\in \cV(k,e)$ for some
$r>0$. For $m\in \N$,
we define
$h_m\in C^\infty_c(\R)$
via
\[
h_m(x):= \frac{r}{m^{k_0}}h(mx).\]
Then $\Supp(h_m)\sub \; ]{-\frac{1}{2m}},\frac{1}{2m}[$
and thus $h_m\in \cV(k,e)$
since, for all $j=0,\ldots, k_0$
and $x\in [-\frac{1}{2},\frac{1}{2}]$,
we have $|h_m^{(j)}(x)|=\frac{rm^j}{m^{k_0}}|h^{(j)}(mx)|<\ve_0$.
We now choose $n\in \N$ such that
$n\geq k_0+2$.
It is easy to construct a
function $\psi\in C^\infty_c(\R)$
such that $\psi(x)=x-n$
for $x$ in some neighbourhood of~$n$ in~$\R$,
and $\Supp(\psi)\sub \;]n-\frac{1}{2}, n+\frac{1}{2}[$.
Then $\phi:=s \cdot \psi\in \cV(k,e)$
for suitable $s>0$.
Choosing $s$ small enough,
we may assume that
$\im(\phi)\sub
[-1,1]$.
The supports of $\phi$ and $h_m$
being disjoint,
we easily deduce from $\phi,h_m\in \cV(k,e)$
that also $\gamma_m:=\phi+ h_m\in \cV(k,e)$.
Then $\gamma_m(0)=0$,
and since $\im(\phi)\sub [-1,1]$,
we have
$f(\gamma_m)(x)
= (h_m\circ \phi)(x)$
for all $x\in W:=\;]n-\frac{1}{2}, n+\frac{1}{2}[$.
For $x\in W$ sufficiently close to~$n$,
we have $\phi(x)=s\cdot (x-n)\in [-\frac{1}{4m},\frac{1}{4m}]$
and thus
$f(\gamma_m)(x)= r\cdot m\cdot s^{k_0+1}\cdot (x-n)^{k_0+1}$,
whence $f(\gamma_m)^{(k_0+1)}(n)=r\cdot m\cdot s^{k_0+1}
\cdot(k_0+1)!\,$.
Thus $f(\gamma_m)\not\in V$ for all
$m\in \N$
such that
$r\cdot m\cdot s^{k_0+1}\cdot (k_0+1)!\geq 1$,
and so $f(\cV(k,e))\!\not\sub \!V$.
As $k$ and $e$ were arbitrary,
(b)\,follows.
\end{proof}
Note that $\Supp(f(\gamma))\sub \Supp(\gamma)$
here, for all $\gamma\in C^\infty_c(\R)$.
\begin{rem}
Although the
map~$f$ from Proposition~\ref{prototype}
is discontinuous and thus not smooth in the
Michal-Bastiani sense,
it is easily seen to be smooth
in the sense of convenient
differential calculus (as any map~$f$
on a ``regular'' countable strict direct
limit $E=\dl\, E_n$\vspace{-1.3 mm} of complete
locally convex spaces,
all of whose restrictions $f|_{E_n}$
are smooth).\footnote{Regularity means that every
bounded subset of~$E$ is contained and bounded
in some~$E_n$.}
\end{rem}
\section{Discontinuous mappings on
{\boldmath $C^\infty_c(M,E)$}}\label{gencase}
In this section,
we generalize our discussion of $C^\infty_c(\R)$
from Section~\ref{secline}
to the spaces $C^\infty_c(M,E)=\dl\, C^\infty_K(M,E)$\vspace{-.8 mm}
of compactly supported
smooth mappings
on a $\sigma$-compact
finite-dimensional
smooth manifold~$M$
with values in
a locally convex space~$E$.
We show:
\begin{prop}\label{propzero}
If $E\not=\{0\}$,
the manifold $M$ is non-compact,
and $\dim(M)>0$,
then there exists a mapping
$f\!: C^\infty_c(M,E)\to C^\infty_c(M,\R)$
such that
\begin{itemize}
\item[\n (a)]
The restriction of $f$ to $C^\infty_K(M,E)$
is smooth, for each compact subset
$K$ of~$M$.
\item[\n (b)]
$f$ is discontinuous at~$0$.
\end{itemize}
In particular,
the locally convex direct limit topology on
$C^\infty_c(M,E)=\dl\, C^\infty_K(M,E)\vspace{-.8 mm}$ is properly
coarser than the topology making
$C^\infty_c(M,E)$ the direct limit of the spaces $C^\infty_K(M,E)$
in the category of topological spaces.
\end{prop}
Instead of proving this
proposition directly, we establish an analogous result
for spaces of sections in bundles of locally
convex spaces, which is no harder to prove.
Noting that the function space $C^\infty_c(M,E)$
is topologically isomorphic to
the space $C^\infty_c(M,M\times E)$
of compactly supported smooth sections
in the trivial bundle $\pr_M\!:M\times E\to M$,
clearly Proposition~\ref{propzero}
is covered by the ensuing discussions for vector
bundles. For background material
concerning bundles of locally convex spaces
and the associated spaces of sections,
the reader is referred to~\cite{SEC}
(or also \cite[Appendix~F]{ZOO}).

For the present purposes, we recall:
if $\pi\!: E\to M$ is a smooth bundle of locally convex
spaces over the finite-dimensional, $\sigma$-compact smooth manifold~$M$,
with typical fibre
the locally convex space~$F$,
then one considers on the space $C^\infty(M,E)$ of all smooth
sections the initial topology with respect
to the family of mappings
\[
\theta_\psi\!: C^\infty(M,E)\to C^\infty(U,F),\;\;\;
\theta_\psi(\sigma):=\sigma_\psi:=\pr_F\circ \psi\circ \sigma|_U^{\pi^{-1}(U)}
\, ,
\]
which take a smooth section~$\sigma$ to its
local representation $\sigma_\psi\!: U\to F$
with respect to the local trivialization
$\psi\!: \pi^{-1}(U)\to U\times F$ of~$E$.
Given a compact subset~$K\sub M$,
the subspace $C^{\,\infty}_K(M,E)\sub C^\infty(M,E)$
of sections vanishing off~$K$ is equipped
with the induced topology,
and $C^\infty_c(M,E):=\bigcup_K C^{\, \infty}_K(M,E)=\dl
\, C^{\,\infty}_K(M,E)$\vspace{-1.3 mm}
is given the locally convex direct limit topology.
\begin{thm}\label{fingeruebung}
Let $M$ be a $\sigma$-compact,
non-compact, finite-dimensional smooth
manifold of dimension $\dim(M)>0$,
and $\pi\!: E\to M$ be a smooth bundle
of locally convex spaces over~$M$,
whose typical fibre is a locally convex
topological vector space~$F\not=\{0\}$.
Then there exists a discontinuous
mapping
$f\!: C^\infty_c(M,E)\to C^\infty_c(M,\R)$
whose restriction to
$C^{\,\infty}_K(M,E)$ is
smooth, for each compact
subset~$K$ of~$M$.
\end{thm}
\begin{proof}
Let $d:=\dim(M)$.
Since~$M$ is non-compact,
there exists a sequence $(U_n)_{n\in \N_0}$
of mutually disjoint
coordinate neighbourhoods $U_n\sub M$
diffeomorphic to $\R^d$
such that local trivializations
$\psi_n\!:\pi^{-1}(U_n)\to U_n\times F$
of~$E$ exist, and such that
every compact subset of~$M$
meets only finitely many of the sets~$U_n$.
We define
\[
\theta_{\psi_n}\!:
C^\infty_c(M,E)\to C^\infty(U_n,F),\;\;\;\;
\theta_{\psi_n}(\sigma):=\sigma_{\psi_n}:=
\pr_F\circ \psi_n\circ \sigma|_{U_n}^{\pi^{-1}(U_n)}\,.
\]
By definition of the topology on $C^\infty_c(M,E)$,
the linear maps $\theta_{\psi_n}$
are continuous.
For each $n\in\N_0$, let
$\kappa_n\!: U_n\to \R^d$ be a $C^\infty$-diffeomorphism;
define $x_n:=\kappa_n^{-1}(0)$.
We choose a function
$h\in C^\infty_c(\R^d,\R)$
such that $h|_{[-1,1]^d}=1$;
we define $h_n\in C^\infty_c(M,\R)$ via
$h_n(x):=h(\kappa_n(x))$ if $x\in U_n$,
$h_n(x):=0$ if $x\in M\,\take\, U_n$.
Let $K_n:=\Supp(h_n)\sub U_n$.
We choose a continuous linear functional
$0\not= \lambda\in F'$,
and pick $v\in F$ such that $\lambda(v)=1$.
Note that $A:=\bigcup_{n\in \N} K_n$
is closed in~$M$, the sequence $(K_n)_{n\in\N}$
of compact sets being locally finite.
Let $\mu\!:\R\times F\to F$ be the
scalar multiplication.
The eventual definition of the mapping~$f$
we are looking for will involve
the map $\Phi\!: E\to M\times \R$, defined via
\begin{equation}\label{dfn1}
\Phi|_{\pi^{-1}(U_n)}:= (\pi|_{\pi^{-1}(U_n)},
\lambda\circ \mu\circ ((h_n\circ \pi)|_{\pi^{-1}(U_n)},
\pr_F\circ \psi_n))
\end{equation}
for $n\in \N$,
and
$\Phi|_{E\take \pi^{-1}(A)}:=(\pi|_{E\take \pi^{-1}(A)}, 0)$.
Note that $\Phi$ is well-defined
as the function in Equation\,(\ref{dfn1})
coincides with $(\pi,0)$
on the set $\bigcup_{n\in \N} \pi^{-1}(U_n \,\take\, A)$.
Also note that $\Phi$ is a fibre-preserving mapping
from~$E$ into the trivial bundle $M\times \R$.
Furthermore,
it is readily verified that~$\Phi$
is a smooth.
By \cite[Thm.\,5.9]{SEC}
(or \cite[Rem.\,F.25\,(a)]{ZOO}), 
the pushforward
\[
C^\infty_c(M,\Phi)\!: C^\infty_c(M,E)\to C^\infty_c(M,M\times \R),\;\;\;\;
\sigma\mto \Phi\circ \sigma\]
is smooth.
For later use,
we introduce the continuous linear map
\[
\Lambda\; :=\; \theta_{\text{id}_{M\times \R}}\!:
C^\infty_c(M,M\times \R)\to C^\infty(M,\R)\, .
\]
Let $\iota\!:\R\to \R^d$
denote the embedding $t\mto (t,0,\ldots, 0)$.
The mapping~$f$ to be constructed will also involve
the map
$\Psi\!: C^\infty_c(M,E)\to C^\infty(\R,\R)$ defined via
\[
\Psi:=
C^\infty(\R,\lambda)
\circ C^\infty(\kappa_0^{-1}\circ \iota ,F)\circ \theta_{\psi_0},
\]
where the pullback
$C^\infty(\kappa_0^{-1}\circ \iota, F)\!:
C^\infty(U_n,F)\to
C^\infty(\R,F)$,
$\gamma\mto \gamma\circ \kappa_0^{-1}\circ \iota$
and the pushforward $C^\infty(\R,\lambda)\!: C^\infty(\R,F)\to C^\infty(\R,
\R)$, $\gamma\mto \lambda\circ \gamma$
are
continuous linear mappings and thus smooth,
by \cite[La.\,3.3, La.\,3.7]{GCX}.
Being a composition of smooth maps,
$\Psi$
is smooth.
We now define the desired map $f\!:C^\infty_c(M,E)\to
C^\infty_c(M,\R)$ via
\[
f\, :=\, \Gamma\circ (\Psi,\Lambda\circ C^\infty_c(M,\Phi)) \;
- \;\lambda\circ \ev_{x_0}\circ \theta_{\psi_0}
\]
(co-restricted from $C^\infty(M,\R)$ to $C^\infty_c(M,\R)$),
where
\[
\Gamma\!: C^\infty(\R,\R)\times C^\infty(M,\R)\to C^\infty(M,\R),\;\;\;\;
\Gamma(\gamma,\eta):=\gamma\circ \eta\]
denotes composition,
and $\ev_{x_0}\!: C^\infty(U_0,F)\to F$
the evaluation map
$\gamma\mto \gamma(x_0)$. Here $\lambda\circ \ev_{x_0}\circ
\theta_{\psi_0}$ is a continuous
linear map and thus smooth.
Explicitly, for $\sigma\in C^\infty_c(M,E)$
\begin{eqnarray*}
f(\sigma)(x) &=&
\Big(\lambda\circ\sigma_{\psi_0}\circ\kappa_0^{-1}\circ\iota\Big)
\bigl( \lambda\bigl( h_n(x)\, \sigma_{\psi_n}(x)\bigr)\bigr)\\
&=&
\lambda\Big( \sigma_{\psi_0} \bigl(\kappa_0^{-1}(
h_n(x)\cdot
\lambda(\sigma_{\psi_n}(x)),\; 0)\bigr)\Big)
\;-\;\lambda(\sigma_{\psi_0}(x_0))
\end{eqnarray*}
if $x\in U_n$ ($n\in\N$),
whereas $f(\sigma)(x)=0$ if $x\in M\,\take\, A$.\vspace{1.3mm}

{\em Claim\/}: {\em
The restriction of~$f$ to $C^\infty_K(M,E)$ is smooth,
for each compact subset~$K$ of~$M$.}\\
To see this, note that
$f(C^{\,\infty}_K(M,E)\sub C^{\,\infty}_K(M,\R)$,
where $C^{\, \infty}_K(M,\R)$ is a closed vector
subspace of $C^\infty(M,\R)$ and
$C^\infty_c(M,\R)$.
Thus, it suffices to show
that $f|_{C^\infty_K(M,E)}$ is smooth
as a map into $C^\infty(M,\R)$
(\cite[Prop.\,1.9]{SEC}, or \cite[La.\,10.1]{BGN}).
But this follows from the Chain Rule,
as $\Gamma$ is smooth by Lemma~\ref{La1}
and also the other constituents of~$f$ are smooth.\vspace{1.3mm}

{\em Claim\/}: {\em $f$ is discontinuous at the zero-section
$\sigma=0$.}
To see this, consider the
set~$V$ of all $\gamma\in C^\infty_c(M,\R)$
such that, for all $n\in \N$
and multi-indices $\alpha\in \N_0^d$
of order $|\alpha|\leq n$,
we have $|\partial^\alpha(\gamma\circ \kappa_n^{-1})(0)|<1$.
It is easily verified that $V$ is a symmetric,
convex zero-neighbourhood
in $C^\infty_c(M,\R)$.
Let $U$ be any convex zero-neighbourhood
in $C^\infty_c(M,E)$;
we claim that $f(U)\not\sub V$.
To see this, set $L_n:=
\kappa_n^{-1}([-1,1]^d)$
for $n\in \N_0$.
Then
\[
\rho_n\!:
C^\infty_{L_n}(M,E)\to C^\infty_{[-1,1]^d}(\R^d,F),\;\;\;\;
\sigma\mto \sigma_{\psi_n}\circ \kappa_n^{-1}
\]
is a topological isomorphism
(cf.\ \cite[La.\,3.9, La.\,3.10]{SEC}
or \cite[La.\,F.9, La.\,F.15]{ZOO})
whose inverse gives rise
to a topological embedding $j_n\!: C^\infty_{[-1,1]^d}(\R^d,F)
\to C^\infty_c(M,E)$.
The linear mapping $\phi\!:\R\to F$, $t\mto tv$
gives rise to a continuous linear map
$C^\infty_{[-1,1]^d}(\R^d,\phi)\!:C^\infty_{[-1,1]^d}(\R^d,\R)
\to C^\infty_{[-1,1]^d}(\R^d,F)$,
$\gamma\mto \phi\circ \gamma$.
Then $W_n:=(j_n\circ C^\infty_{[-1,1]^d}(\R^d,\phi))^{-1}(\frac{1}{2}U)$
is a convex zero-neighbourhood
in $C^\infty_{[-1,1]^d}(\R^d,\R)$.
Thus, there exists $k_n\in \N_0$
and $\ve_n>0$ such that $W_{k_n,\ve_n}\sub W_n$,
where $W_{k_n,\ve_n}$
is the set of
all $\gamma\in C^\infty_{[-1,1]^d}(\R^d,\R)$
such that $\sup\{|\partial^\alpha \gamma(x)|\!:
x\in [-1,1]^d\}<\ve_n$
for all $\alpha\in \N_0^d$
such that $|\alpha|\leq k_n$.
We let $g\in C^\infty_{[-1,1]^d}(\R^d,\R)$
be a function such that
$g(y_1,\ldots,y_d)=y_1^{{k_0}+1}$
for all $y=(y_1,\ldots, y_d)\in [-\frac{1}{2},\frac{1}{2}]^d$.
Then $rg\in W_{k_0,\ve_0}$ for some $r>0$.
It is clear from the definition of $W_{k_0,\ve_0}$
that then also $\gamma_m\in W_{k_0,\ve_0}$
for all $m\in \N$,
where
\[
\gamma_m\!:\R^d\to\R\, ,\quad
\gamma_m(y_1,\ldots,y_d):= \frac{r}{m^{k_0}}\,g(my_1,y_2,\ldots, y_d)\, .
\]
Thus $\tau_m:= j_0(\phi\circ \gamma_m)\in \frac{1}{2}U$.

Let $\ell:=k_0+1$;
we easily find $\eta\in W_{k_\ell,\ve_\ell}$
such that, for suitable $s>0$, we have
$\eta(y)=s\cdot y_1$ for $y=(y_1,\ldots, y_d)$
in some zero-neighbourhood in~$\R^d$. We
define $\tau:=j_\ell(\phi\circ \eta)\in \frac{1}{2}U$.
Then $\sigma_m:=\tau_m+\tau\in U$
by convexity of~$U$.
Consider $g_m:=f(\sigma_m)\circ \kappa_\ell^{-1}\!:
\R^d\to \R$.
For $y\in [-1,1]^d$
sufficiently close to~$0$, we have $\eta(y)=sy_1$
and $m|\eta(y)|\leq\frac{1}{2}$. Thus
\[
g_m(y)=\gamma_m(\eta(y),0,\ldots,0)
=r\cdot m\cdot s^{k_0+1}\cdot y_1^{k_0+1},
\]
entailing
that
$\frac{\partial^{k_0+1} g_m}{\partial y_1^{k_0+1}}(0)
=r\cdot m\cdot s^{k_0+1}\cdot (k_0+1)!\,$.
Hence $f(\sigma_m)\not\in V$
for each $m\in \N$
such that
$r\cdot m\cdot s^{k_0+1}\cdot (k_0+1)!\geq 1$.
We have shown that $f(U)\not\sub V$
for any $0$-neighbourhood~$U$
in $C^\infty_c(M,E)$, although
$f(0)=0$. Thus $f$ is discontinuous
at $\sigma=0$.
\end{proof}
\section{Further examples}\label{secmisc}
We describe various pathological bilinear mappings.
\begin{prop}
Let $\K\in\{\R,\C\}$. The pointwise multiplication map
\[
\mu\!:C^\infty(\R,\K)\times C^\infty_c(\R,\K)\to C^\infty_c(\R,\K),\;\;\;\;
\mu(\gamma,\eta):=\gamma\cdot\eta
\]
is a hypocontinuous bilinear $($and thus sequentially
continuous$)$ mapping on the locally convex direct limit
\[
C^\infty(\R,\K)\times C^\infty_c(\R,\K)=
\dl\, (C^\infty(\R,\K)\times C^\infty_{[-n,n]}(\R,\K))\, ,\]
whose restriction to $C^\infty(\R,\K)\times C^\infty_{[-n,n]}(\R,\K)$
is continuous bilinear and thus $\K$-analytic,
for each $n\in \N$.
However, $\mu$ is discontinuous.
\end{prop}
\begin{proof}
Using the Leibniz Rule for the differentiation
of products of functions,
it is easily verified that~$\mu$ is
separately continuous.\footnote{Alternatively,
we can obtain the assertion as a special
case of \cite[Cor.\,2.7]{SEC}
or \cite[La.\,4.5\,(a) and Prop.\,4.19\,(d)]{ZOO},
combined with the locally convex direct limit property.}
The spaces $C^\infty(\R,\K)$ and $C^\infty_c(\R,\K)$
being barrelled,
this entails that~$\mu$ is hypocontinuous
and thus sequentially continuous
\cite[Thm.\,41.2]{Tre}.
The restriction of~$\mu$ to
$C^\infty(\R,\K)\times C^\infty_{[-n,n]}(\R,\K)$
is a sequentially continuous
bilinear mapping on a product of metrizable
spaces and therefore continuous.
To see that $\mu$ is discontinuous,
consider the zero-neighbourhood
\[
W:=\{\gamma\in C^\infty_c(\R,\K)\!: \; (\forall x\in \R)\;\,|\gamma(x)|<1\}
\]
in $C^\infty_c(\R,\K)$.
If $U$ is any zero-neighbourhood
in $C^\infty(\R,\K)$
and $V$ any zero-neighbourhood in $C^\infty_c(\R,\K)$,
then there exists a compact subset
$K$ of~$\R$ such that
\[
(
\forall \gamma\in C^\infty(\R,\K))\;\;\;
\gamma|_K=0\; \impl\; \gamma\in U.
\]
Pick any $x_0\in \R\, \take\, K$.
There is a function $\phi\in C^\infty_c(\R,\K)$
such that $\phi(x_0)\not=0$ and $\Supp(\phi)\sub
\R\,\take\, K$.
Then $r\phi\in V$ for some $r>0$,
and $t\phi\in U$ for all $t\in \R$.
Choosing $t\geq  \frac{1}{r\cdot|\phi(x_0)|^2}$,
we have $(t\phi,r\phi)\in U\times V$ but
$|\mu(r\phi,t\phi)(x_0)|=rt|\phi(x_0)|^2\geq 1$,
entailing that $\mu(U\times V)\not\sub W$.
Thus $\mu$ is discontinuous at~$(0,0)$.
\end{proof}
Another instructive example is the following
(compare also the examples in \cite{DaW}):
\begin{example}
Let $E_1\subset  E_2\subset \cdots$ be a strictly ascending
sequence of Banach spaces, such that $E_{n+1}$
induces the given topology on~$E_n$.
Set $E:=\dl\, E_n$\vspace{-.8mm}
and $F:=E'_b$. For example,
we can take $E_n:=L^2[-n,n]$,
in which case $E=L^2_\comp(\R)$ and
$F=L^2_\loc(\R)=\pl\, L^2[-n,n]$.
Then
$A_n:=F\times E_n\times \K\times \K$
is a Fr\'{e}chet space (and reflexive
in the example $E_n=L^2[-n,n]$).
The evaluation map $E_n'\times E_n\to \R$ being continuous
as~$E_n$ is a Banach space,
it is easy to see that~$A_n$
becomes a unital associative
topological algebra via
\begin{equation}\label{formull}
(\lambda_1,x_1,z_1,c_1)\cdot(\lambda_2,x_2,z_2,c_2)
\,:=\,
\bigl(
c_1\lambda_2+c_2\lambda_1,\, c_1x_2+c_2x_1,\,c_1z_2+ \lambda_1(x_2)+z_1c_2,
\, c_1c_2 \bigr)\, .
\end{equation}
The
multiplication can be visualized by considering
$(\lambda,x,z,c)\!\in \!A_n$ as the 3-by-3 matrix
{\scriptsize
\[
\left(
\begin{array}{ccc}
c & \lambda & z\\
0 & c & x \\
0 & 0 & c
\end{array}
\right).
\]
}The topological algebras~$A_n$
are very well-behaved:
they have open groups of units,
and inversion is a $\K$-analytic map.
We can also use Formula\,(\ref{formull})
to define a multiplication map
$\mu\!: A\times A\to A$
turning the direct limit locally convex space
$A:=F\times E\times \K\times \K=\dl\, A_n$\vspace{-.8 mm}
into a unital, associative algebra.
However, although
the restriction of~$\mu$
to $A_n\times A_n$ is a continuous bilinear map
for each~$n\in\N$,
$\mu\!: A\times A=\dl\, (A_n\times A_n)\to A$\vspace{-1.3 mm}
is discontinuous (since the evaluation map
$E'_b\times E\to\R$ is discontinuous,
the space~$E$ not being normable).
We refer to \cite[Section~10]{Glo}
for more details.
\end{example}
\appendix
\section*{Appendix: Proof of Lemma~\ref{La1}}
We give a proof which is
as elementary as possible,
by reducing the assertion to
the case $M=\R^d$.
First, let $M$ be a finite-dimensional,
$\sigma$-compact smooth manifold,
of dimension~$d$.
We choose
an open cover $(U_j)_{j\in J}$ of~$M$
and $C^\infty$-diffeomorphisms
$\kappa_j\!: U_j\to\R^d$.
Then
\[
\Phi : \, C^\infty(M,\R^m)\, \to \, \prod_{j\in J} C^\infty(\R^d,\R^m)\,
=:\, P\,,\qquad 
\Phi(\gamma)\, :=\, (\gamma\circ \kappa_j^{-1})_{j\in J}
\]
is a topological embedding onto a closed
vector subspace of the cartesian product~$P$
(cf.\ \cite[La.\,3.7]{SEC}).
Therefore $\Gamma$ is smooth
if and only if $\Phi\circ \Gamma$ is smooth
(\cite[Prop.\,1.9]{SEC} or \cite[La.\,10.1]{BGN}),
if and only if each component $\pr_j\circ \Phi\circ \Gamma$
is smooth \cite[La.\,10.3]{BGN}, where $\pr_j\!: P\to C^\infty(\R^d,\R^m)$
is the projection onto the $j$-coordinate.
But
\[
\pr_j(\Phi(\Gamma))(\gamma,\eta)\; =\;
\gamma\circ \eta\circ \kappa_j^{-1}
\; =\;\tilde{\Gamma}\bigl(\gamma, C^\infty(\kappa_j^{-1}, \R^n)(\eta)\bigr)
\]
for all $\gamma\in C^\infty(\R^n,\R^m)$ and $\eta\in C^\infty(M,\R^n)$,
where
\[
\tilde{\Gamma}\!: C^\infty(\R^n,\R^m)\times C^\infty(\R^d,\R^n)\to
C^\infty(\R^d,\R^m)
\]
is the composition map
and $C^\infty(\kappa_j^{-1},\R^n)\!:
C^\infty(M,\R^n)\to C^\infty(\R^d,\R^n)$, $\eta\mto\eta\circ\kappa_j^{-1}$
is continuous linear and thus smooth, by
\cite[La.\,3.7]{GCX}.
Hence $\pr_j\circ \Phi\circ \Gamma$
(and thus $\Gamma$) will be smooth if so is~$\tilde{\Gamma}$.\\[3mm]
By the reduction step just performed,
it only remains to prove Lemma~\ref{La1}
for $M=\R^d$, which we assume now.
We show by induction
on $k\in \N_0$
that~$\Gamma$ is~$C^k$.\\[3mm]
{\em The case $k=0$}.
Let $\gamma\in C^\infty(\R^n,\R^m)$,
$\eta \in C^\infty(\R^d,\R^n)$
and $(\gamma_i,\eta_i)_{i\in \N}$
be a sequence in $C^\infty(\R^n,\R^m)\times
C^\infty(\R^d,\R^n)$
converging to $(\gamma,\eta)$.
We have
to show that
\begin{eqnarray}
\delta_i & := & \Gamma(\gamma_i,\eta_i)-
\Gamma(\gamma,\eta)\; =\;
\gamma_i\circ \eta_i
-\gamma\circ \eta\nonumber \\
& = &
(\gamma_i-\gamma)\circ \eta_i
\; + \; (\gamma\circ
\eta_i-\gamma\circ \eta)\label{lastline}
\end{eqnarray}
converges to~$0$ in $C^\infty(\R^d,\R^m)$.
To see this, we first check convergence
in $C^0(\R^d, \R^m)$
(equipped with the topology of uniform convergence
on compact sets).
Given a compact set $K\sub \R^d$,
the set
$\bigcup_{i\in \N}\eta_i(K)$
is bounded and
hence has compact closure~$L$ in~$\R^d$.
Now the
first term in~(\ref{lastline})
converges uniformly to~$0$
on~$K$
since $\gamma_i-\gamma\to 0$
uniformly on~$L$ as $i\to\infty$.
The
second term converges
uniformly to~$0$ on~$K$
since $\gamma|_L$ is uniformly continuous
and $\eta_i \to \eta$ uniformly on~$K$.
Using the Chain Rule,
for each fixed
multi-index $\alpha\in \N_0^d$
of order $\geq 1$,
we find
polynomials
$P_\beta \in \R[(X_\gamma)_{\gamma\leq \alpha}]$
in indeterminates~$X_\gamma$,
for multi-indices $\beta\in \N_0^n$
of order $|\beta| \leq |\alpha|$,
such that
\begin{eqnarray*}
\partial^\alpha
\delta_i& = &
\sum_{|\beta|\leq |\alpha|}
((\partial^\beta\gamma_i-\partial^\beta\gamma)\circ \eta_i)\cdot
P_\beta((\partial^\gamma\eta_i)_{\gamma\leq \alpha})\\
& &
+\sum_{|\beta|\leq |\alpha|}
(\partial^\beta\gamma \circ \eta_i)\cdot
(P_\beta((\partial^\gamma\eta_i)_{\gamma\leq \alpha})-
P_\beta((\partial^\gamma\eta)_{\gamma\leq \alpha}))\\
& &
+\sum_{|\beta|\leq |\alpha|}
(\partial^\beta\gamma\circ \eta_i-
\partial^\beta\gamma\circ \eta)\cdot
P_\beta((\partial^\gamma\eta)_{\gamma\leq\alpha}).
\end{eqnarray*}
We easily deduce from this formula
that $\partial^\alpha\delta_i$ converges
to~$0$ as $i\to\infty$, uniformly
on compact sets.
We have shown that $\delta_i\to 0$ in
$C^\infty(\R^d,\R^m)$.
Thus $\Gamma$ is continuous.\vspace{2 mm}\\
{\em Induction step}.
Suppose that $\Gamma$ is of class~$C^k$,
where $k\in\N_0$.
Given $\gamma,\gamma_1\in C^\infty(\R^n,\R^m)$,
$\eta,\eta_1\in C^\infty(\R^d,\R^n)$,
we have
\begin{equation}\label{kk}
{\textstyle
\frac{1}{t}\left(\Gamma(\gamma+t\gamma_1,\eta+t\eta_1)
-\Gamma(\gamma,\eta)\right)
\; = \;
\frac{1}{t}\left(\gamma\circ (\eta+t\eta_1)-\gamma\circ\eta\right)
\;+\;\gamma_1\circ (\eta+t\eta_1)}
\end{equation}
for $0\not=t\in \R$.
Given $t\in \R$, define $F_t\!:\R^d\to \R^m$ via
\[
F_t(x)\; :=\; \int_0^1H(x,st)\; ds\, ,
\]
where $H\!: \R^d \times \R\to \R^m$,
$H(x,r):=
d\gamma(\eta(x)+r\eta_1(x);\eta_1(x))$.
Clearly~$H$
is smooth.
It is easy to see
that $F_t(x)\to F_0(x)$
uniformly for $x$ in a compact set, as $t\to 0$.
Furthermore, differentiating under
the integral sign we find
that
$\partial^\alpha F_t(x)=\int_0^1 \partial^{(\alpha,0)}H(x,st)\,ds$
for $\alpha\in \N_0^d$,
which converges uniformly for $x$ in a compact set
to $\partial^\alpha F_0(x)$ as $t\to 0$.
Since
\[
F_t\; =\;
{\textstyle
\frac{1}{t}\left(\gamma\circ (\eta+t\eta_1)-\gamma\circ\eta\right)}
\]
for $t\not=0$, by the Mean Value Theorem,
we see
that the first term on the right hand side
of (\ref{kk})
converges
to $F_0= (d\gamma)\circ (\eta,\eta_1)
= \wt{\Gamma}(d\gamma,(\eta,\eta_1))$
in $C^\infty(\R^d,\R^m)$ as $t\to 0$,
where
$\wt{\Gamma}\!: C^\infty(\R^n\times \R^n,\R^m)\times C^\infty
(\R^d,\R^n\times \R^n)\to C^\infty(\R^d,\R^m)$ is the composition map.

To tackle the second term,
define
$G_t:=\gamma_1\circ (\eta+t\eta_1)=\Gamma(\gamma_1,\eta+t\eta_1)$
for $t\in\R$.
Since $\Gamma$ is continuous
by the above, we have $G_t\to G_0=\gamma_1\circ \eta$
in $C^\infty(\R^d,\R^m)$ as $t\to 0$.
Thus the second term in Equation\,(\ref{kk})
converges to~$\gamma_1\circ \eta$.

Summing up,
we have shown that
$d \Gamma(\gamma,\eta;\gamma_1,\eta_1)$ exists,
and is given by
\begin{equation}\label{facilit}
d\Gamma(\gamma,\eta;\gamma_1,\eta_1)
\; =\; \wt{\Gamma}(d\gamma,(\eta,\eta_1)) + \Gamma(\gamma_1, \eta).
\end{equation}
The
map $C^\infty(\R^n,\R^m)\to C^\infty(\R^n\times\R^n,\R^m)$,
$\gamma\mto d \gamma$ is continuous linear
(cf.\ \cite[La.\,3.8]{GCX}),
and $\Gamma$, $\wt{\Gamma}$ are~$C^k$, by induction.
Hence Equation\,(\ref{facilit}) shows that
$d\Gamma$ is~$C^k$.
Thus $\Gamma$ is~$C^{k+1}$,
as required.\\[3mm]
{\em Acknowledgement.} The research was partially supported by
DFG, FOR~363/1-1.

{\small Helge Gl\"{o}ckner, TU Darmstadt, Fachbereich Mathematik AG~5,
Schlossgartenstr.\,7,\\
64289 Darmstadt, Germany. E-Mail: gloeckner\at{}mathematik.tu-darmstadt.de}
\end{document}